\documentclass[a4paper, 12pt, reqno]{amsart}
\usepackage{amsmath,amsthm,amssymb,amsfonts, bm, url, a4wide}
\usepackage{hyperref}
\newcommand{\R}{\mathbb{R}}
\newcommand{\Rp}{{\R}^{+}}    

\newcommand{\Z}{\mathbb{Z}}
\newcommand{\cH}{\mathcal{H}}
\newcommand{\N}{\mathbb{N}}
\newcommand{\Q}{\mathbb{Q}}
\newcommand{\cA}{\mathcal{A}}

\newcommand{\vv}[1]{{\mathbf{#1}}}
\def\ds{\displaystyle}

\begin {document}
\newtheorem{pretheorem}{Theorem}

\newtheorem{theorem}{Theorem}
\newtheorem{lemma}{Lemma}
\newtheorem*{KJ}{Theorem (Khintchine--Jarn\'{\i}k)}
\newtheorem*{theoremBVV}{Theorem BVV}

\newtheorem*{theoremG}{Theorem G}
\newtheorem{corollary}{Corollary}

\title{On weighted inhomogeneous Diophantine approximation on planar curves}

\author[M.~Hussain]{Mumtaz Hussain}

\address{Mumtaz Hussain, Department of Mathematics and Statistics, La Trobe University, Melbourne, 3086, Victoria, Australia}
\email{m.hussain@latrobe.edu.au}

\author{Tatiana Yusupova}

\address {Tatiana Yusupova, Department of Mathematics, University of York, Heslington, York,
YO105DD, UK}
\email {tatiana.yusupova@gmail.com}


\begin{abstract}
This paper develops the metric theory of simultaneous inhomogeneous Diophantine approximation on a planar curve with respect to multiple approximating functions. Our results naturally  generalize  the homogeneous Lebesgue measure and Hausdorff dimension results for the sets of simultaneously well-approximable points on planar curves, established  in Badziahin and Levesley  ({\em Glasg. Math. J.}, 49(2):367--375, 2007), Beresnevich et al.  ({\em Ann. of Math. (2)}, 166(2):367--426, 2007), Beresnevich and Velani  ({\em Math. Ann.}, 337(4):769--796, 2007) and Vaughan and Velani ({\em Invent. Math.}, 166(1):103--124, 2006). 
\end{abstract}
\maketitle
\section{Introduction}

Throughout $\psi_1, \psi_2:\N\to\R^+$ are  monotonic functions  such that $\psi_i(q)\to 0$ as $q\to\infty$ for each $i$ and will be referred to as {\em approximating functions}. When $\psi_1=\psi_2$ we will denote them as $\psi.$

 Let $\bm\theta=(\theta_1,\theta_2)\in\R^2$ be fixed throughout and $\mathcal{A}(\psi_1,\psi_2,\bm\theta)$ denote the set of \linebreak $\vv x=(x_1,x_2)\in\R^2$ such that the system of inequalities
$$
\begin{cases}
 ||q x_{1}-\theta_1|| < \psi_1(q) \\[2ex]
 ||q x_{2}-\theta_2|| < \psi_2(q)
\end{cases}
$$
holds for infinitely many $q\in\N$. Here and throughout $||\cdot||$ denotes the distance to the nearest integer.  The points in $\mathcal{A}(\psi_1,\psi_2,\bm\theta)$ are said to be {\em simultaneously $(\psi_1,\psi_2)$-approximable in the inhomogeneous setup}. When $\psi_1=\psi_2$ the set $\mathcal{A}(\psi_1,\psi_2,\bm\theta)$ is denoted by $\mathcal{A}(\psi,\bm\theta)$. In the case of the inhomogeneous factor $\bm\theta$ being equal to $\bm 0$, the corresponding set $\cA(\psi)$ will be the usual homogeneous set of simultaneously $\psi$-approximable points on the plane.
\smallskip

The following statement provides an elegant yet simple criterion for the ‘size’ of $\cA(\psi,\bm\theta)$ expressed in terms of $s$-dimensional Hausdorff measure $\cH^s$. It is a generalization of the classical theorems of Khintchine (\cite{khintchine-first}, 1924) and Jarn\'{\i}k (\cite{Jarnik}, 1931). For basic concepts and definitions  we refer to \cite{berdod, Cassel, Falc}.

\begin{KJ}
Let $\psi:\N\to\R^+$ be a monotonic approximating function, $\bm\theta\in\R^2$ and $s\in(0,2]$. Then
$$\cH^s(\mathcal{A}(\psi,\bm\theta))=\begin {cases}
 0 \ & {\rm if } \quad \sum\limits_{q=1}^{\infty}q^{2-s}\psi^{s}(q)<\infty \\[2ex]
 \cH^s(\R^2) \ & {\rm if } \quad \sum\limits_{q=1}^{\infty}q^{2-s}\psi^{s}(q)=\infty.
\end {cases}$$
\end{KJ}

When $s = 2$, the measure $\cH^2$ is equivalent to two-dimensional Lebesgue measure on the plane and the theorem corresponds to Khintchine's theorem. When $s < 2$, the homogeneous case of the theorem corresponds to Jarn\'{\i}k's theorem and can be regarded as a Hausdorff measure version of Khintchine's theorem. For further details see \cite{BDV_mtl} and references within.

We would like to be able to describe the measure theoretic structure of the set $\mathcal{A}(\psi_1,\psi_2,\bm\theta)$ taken in the intersection with a planar curve $\mathcal{C}$, i.e. the set $\mathcal{C}\cap\mathcal{A}(\psi_1,\psi_2,\bm\theta)$. Note that the dependency between the points of the planar curve introduces major difficulties in attempting to describe this set. Until now only the partial results were obtained. 

Without loss of generality, we will assume that $\mathcal C := \mathcal{C}_f=\{(x, f(x)): \; x\in I\}$ is given as
the graph of a function $f : I\to \R,$ where $I$ is some interval of $\R$. As usual, $C^{(n)}(I)$ will
denote the set of $n$-times continuously differentiable functions defined on some interval $I$
of $\R$. A planar curve $\mathcal C$ is non-degenerate if the set of points on $\mathcal C$ at which the curvature vanishes is a set
of one-dimensional Lebesgue measure zero. The geometric interpretation of non-degeneracy
is that the curve is curved enough that it deviates from any hyperplane.  For further details and consequences of this definition we refer to \cite {VB, BZ, KM}.

A combination of results established in \cite{BDV_main} and \cite{VV} provides a complete analogue of the  Khintchine--Jarn\'{\i}k theorem for a set $\mathcal{C}_f\cap \mathcal{A}(\psi)$. For weighted sets $\mathcal{C}_f\cap \mathcal{A}(\psi_1, \psi_2)$ within the homogeneous setup the results given in \cite{BL} and \cite{BV} provide only a Lebesgue measure analogue of the Khintchine--Jarn\'{i}k theorem.

\begin{pretheorem}[BV 2007 and BL 2007]\label{thB}
Let $\psi_1$ and $\psi_2$ be approximating functions, $f\in C^{(3)}(I)$ and let $\mathcal{C}_f$ be a non-degenerate curve. Then
\[|\mathcal{C}_f\cap \mathcal{A}(\psi_1, \psi_2)|_{\mathcal{C}_f}=
\begin {cases}
 \text{\upshape Zero} \ &{\rm if } \quad\sum\psi_1(q)\psi_2(q) < \infty \\[3ex]
 \text{\upshape  Full} \ &{\rm if } \quad\sum\psi_1(q)\psi_2(q) = \infty,
\end {cases}\]
where $|\ \cdot \ |_{\mathcal{C}_f}$ denotes the induced Lebesgue measure on ${\mathcal{C}_f}$.
\end{pretheorem}

The Hausdorff measure theory for weighted sets near curves remained undeveloped even in this homogeneous case with some Hausdorff dimension results obtained in \cite{BV}.

The main result of \cite{BVV} gives a stronger inhomogeneous Hausdorff measure statement, however, it is restricted to the situation when $\psi_1=\psi_2$.

\begin{pretheorem}[BVV 2010]\label{thA}
Let $\psi:\N\to\R^+$ be an  approximating function, $\bm\theta\in\R^2$ and $f\in C^{(3)}(I)$. Furthermore, let $1/2<s\leqslant 1$ and assume that $\cH^s\{x\in I: \;f''(x)=0\}=0$. Then
$$\cH^s(\mathcal{C}_f\cap \mathcal{A}(\psi,\bm\theta))=\begin {cases}
 0 \ &{\rm if } \quad \sum\limits_{q=1}^{\infty}q^{1-s}\psi^{s+1}(q)<\infty \\[3ex]
 \cH^s(I) \ &{\rm if } \quad \sum\limits_{q=1}^{\infty}q^{1-s}\psi^{s+1}(q)=\infty.
\end {cases}$$
\end{pretheorem}

This shows that nothing is known when $\psi_1\neq\psi_2$ and $\bm\theta\neq (0,0)$.

\section{Statements of the results}

The aim of this paper is to develop the Hausdorff measure theory for $\mathcal{C}_f\cap \mathcal{A}(\psi_1, \psi_2,\bm\theta)$. The dimension results of \cite{BV} will be a consequence of our main result. In order to make our results even more general we will work with the general dimension function $h$ for Hausdorff measure instead of the specific case $h(r)=r^s$ appearing in the previous results.

\noindent Before stating new results, we introduce a useful notion concerning dimension functions.

A \textit{dimension function} is an increasing, continuous function $h:\R^+\rightarrow \R^+$ such that $h(r)\to 0$ as $r\to 0$. It is said to be \textit{regular} if there exist constants $r_0,\lambda_1,\lambda_2\in (0,1)$ such that $h(\lambda_1r)\leqslant\lambda_2h(r)$ for $r\in(0,r_0)$.

\noindent Our main result can now be formulated as follows.

\begin{theorem}\label{t04}
Let $\psi_1$ and $\psi_2$ be approximating functions such that $\max\{\psi_1(q),\psi_2(q)\}\geqslant q^{-\eta}$ for some $\eta <1$. Let $\bm\theta\in\R^2$ and  $f\in C^{(3)}(I)$. Let $h$ be a regular dimension function such that $r^{-1}h(r)\to \infty$ as $r\to 0$ and $r^{-1}h(r)$ is decreasing. Furthermore, assume that $\cH^h(\{x\in I: f''(x)=0\})=0$. Then
\[\cH^h(\mathcal{C}_f \cap \mathcal{A}(\psi_1, \psi_2, \bm\theta))=
\begin {cases}
 0 \ &{\rm if } \quad\sum\limits_{q=1}^{\infty} q\cdot h\left(\frac{\min\{\psi_1(q),\psi_2(q)\}}{q}\right)\cdot\max\{\psi_1(q),\psi_2(q)\} < \infty \\[2ex]
 \infty \ & {\rm if } \quad\sum\limits_{q=1}^{\infty} q\cdot h\left(\frac{\min\{\psi_1(q),\psi_2(q)\}}{q}\right)\cdot\max\{\psi_1(q),\psi_2(q)\} = \infty.
\end {cases}\]

\end{theorem}

\noindent\textbf{Remark 1}: Note that if the dimension function $h$ satisfies the conditions of the theorem, then $\cH^h$ is not comparable to one-dimensional Lebesgue measure. The point is that the dimension function $h(r)=r$ does not satisfy the condition $r^{-1}h(r)\to \infty$ as $r\to 0$. We will address the case of Lebesgue measure separately - see Theorem~\ref{t02}.

\noindent\textbf{Remark 2}: Theorem~\ref{t04} is a natural generalization of Theorem~\ref{thA}. Indeed, with $s$-dimensional Hausdorff measure $\cH^s$ and $\psi_1=\psi_2$, the condition $\eta<1$ imposed in our theorem implies that $s>1/2$.

As a corollary of Theorem~\ref{t04}, the following inhomogeneous generalization of the Hausdorff dimension result from \cite{BV} can be obtained.

\begin{corollary}\label{c05}
For $i=1,2$ let $\psi_i(q)=q^{-v_i}$, where $0<\min\{v_1,v_2\}<1$. Let $\bm\theta=(\theta_1,\theta_2)\in\R^2$ and $f\in C^{(3)}(I)$. Let $$s_0=\frac{2-\min\{v_1,v_2\}}{1+\max\{v_1,v_2\}}$$ and assume that $s_0<1$ and $\cH^{s_0}(\{x\in I: f''(x)=0\})=0$.
Then $$\dim(\mathcal{C}_f \cap \mathcal{A}(\psi_1, \psi_2, \bm\theta))=s_0$$
and $\cH^{s_0}(\mathcal{C}_f \cap \mathcal{A}(\psi_1, \psi_2,\bm\theta))=\infty$.
\end{corollary}

\noindent\textbf{Remark 3}: In the homogeneous case, the dimension part of Corollary~\ref{c05} essentially corresponds to Theorem 4 in \cite{BV}. Indeed, the condition $s_0<1$ is guaranteed by assuming $v_1+v_2> 1$ and $\cH^{s_0}(\{x\in I: f''(x)=0\})=0$ implies that $\dim(\{x\in I: f''(x)=0\})\leqslant s_0$.
The Hausdorff measure part of Corollary~\ref{c05} is new even in the homogeneous setup.

To prove Theorem~\ref{t04} we shall make use of ubiquity and its inhomogeneous analogue as described in \cite{BVV}.
The ideas and techniques used for the proof will also enable us to prove the Lebesgue measure version of Theorem~\ref{t04}. Recall that Lebesgue measure is excluded in Theorem~\ref{t04} - see Remark 1.

\begin{theorem}\label{t02}
Let $\psi_1$ and $\psi_2$ be approximating functions, $\bm\theta\in\R^2$, $f\in C^{(3)}(I)$ and $\mathcal{C}_f$ be a non-degenerate curve. Then
\[|\mathcal{C}_f\cap \mathcal{A}(\psi_1, \psi_2,\bm\theta)|_{\mathcal{C}_f}=
\begin {cases}
 \text{\upshape Zero} \ & {\rm if } \quad \sum\limits_{q=1}^{\infty}\psi_1(q)\psi_2(q) < \infty \\[3ex]
 \text{\upshape Full} \ & {\rm if } \quad \sum\limits_{q=1}^{\infty}\psi_1(q)\psi_2(q) = \infty.
\end {cases}\]
\end{theorem}

\emph{Notation.} To simplify notation in the proofs below the Vinogradov symbols $\ll$ and $\gg$ will be used to indicate an inequality with an unspecified positive multiplicative constant.
If $a\ll b$ and $a\gg b$  we write $a\asymp
b$, and say that the quantities $a$ and $b$ are comparable.

\section{Auxiliary results: ubiquity}\label{ubiquity}
The proofs of the divergent parts of Theorems~\ref{t04} and \ref{t02} are based on proving that the described systems are ubiquitous. The idea and the concept of ubiquity was originally formulated by Dodson, Rynne and Vickers in \cite{DRV} and coincided in part with the concept of regular systems of Baker and Schmidt (\cite{BS}). Both have proven to be extremely useful in obtaining lower bounds for the Hausdorff dimension of limsup sets. The recent further development of ubiquity by Beresnevich, Dickinson and Velani in \cite{BDV_mtl} provides a very general and abstract approach for establishing the Hausdorff measure of a large class of limsup sets. In this section, for the sake of simplicity, we introduce a restricted form of ubiquity to deal with Diophantine approximation on planar curves. We refer to \cite{BDV_main} for detailed acquisition of the ubiquity framework and only present a restricted version of it for our use. Note that so far the ubiquity framework has been used for one approximating function. However, the framework is broad ranging and can be used with respect to multiple approximating functions. We make appropriate changes below.

\subsection{Ubiquitous systems on $\R$}
Let $I$ be an interval in $\R$ and $\mathcal{R}=(R_{\alpha})_{\alpha\in\mathcal{J}}$ be a family of resonant points $R_{\alpha}$ of $I$ indexed by an infinite, countable set $\mathcal{J}$. Here, for $\alpha\in\mathcal{J}$, the resonant points are defined as  $\{x\in I:|x-R_{\alpha}|=0 \}.$ Next let $\beta:\mathcal{J}\to\R^{+}:\alpha\to\beta_{\alpha}$ be a positive function on $\mathcal{J}$. Thus, the function $\beta$ attaches a `weight' $\beta_{\alpha}$ to the resonant point $R_{\alpha}$. Also, for $t\in\N$ let $$J(t)=\{\alpha\in\mathcal{J}:\beta_{\alpha}\leqslant k^t\}$$ for a fixed constant $k>1$. Assume that for every $t$ the set $J(t)$ is finite.

Throughout, $\rho:\R^{+}\to\R^{+}$ will denote a function with $\lim_{t\to\infty}\rho(t)=0$. The function $\rho$ is usually referred to as the ubiquitous function. Also $B(x,r)$ will denote the ball (or rather the interval) centered at $x$ of radius $r$.

\noindent\textsc{\textbf{Definition:} }\textit{
Suppose there exists a function $\rho$ and an absolute constant $\kappa>0$ such that for any interval $J\subseteq I$
$$\lim_{t\to\infty}\inf\big|\bigcup\limits_{\alpha\in J(t)}(B(R_{\alpha},\rho(k^t))\cap J)\big|\geqslant \kappa |J|,$$
then the system $(\mathcal{R},\beta)$ is called \textbf{locally ubiquitous} in $I$ \textbf{with respect to} $\rho$.}

As a consequence of this definition there is a series of measure theoretic results established in \cite {BDV_mtl, BDV_main}. We will only be exploiting the above definition in the framework below.

\subsection{Ubiquitous systems near curves.}
With $n\geqslant 2$, let $\mathcal{R}=(R_{\alpha})_{\alpha\in\mathcal{J}}$ be a family of resonant points $R_{\alpha}$ of $\R^n$ indexed by an infinite set $\mathcal{J}$. Let $\beta:\mathcal{J}\to\R^{+}:\alpha\to\beta_{\alpha}$ is a positive function on $\mathcal{J}$. For a point $R_{\alpha}$ in $\mathcal{R}$ let $R_{\alpha,k}$ represent the $k^{\text{th}}$ coordinate of $R_{\alpha}$. Thus, $R_{\alpha}=(R_{\alpha,1}, R_{\alpha,2},\ldots, R_{\alpha,n})$. We will use the notation $\mathcal{R_C}(\Phi)$ to denote the sub-family of resonant points $R_{\alpha}$ in $\mathcal{R}$ which are ``$\Phi$-close'' to the curve $\mathcal{C}:=\mathcal{C}_f=\{(x,f_2(x),\ldots,f_n(x)):x\in I\}$ where $\Phi$ is an approximating function, $f=(f_1,\ldots,f_n):I\to\R^n$ is a continuous map with $f_1(x)=x$ and $I$ is an interval in $\R$. Formally and more precisely $$\mathcal{R_C}(\Phi)=(R_{\alpha})_{\alpha\in\mathcal{J_C}(\Phi)},$$ where
$$\mathcal{J_C}(\Phi)=\{\alpha\in\mathcal{J}:\max_{1\leqslant k\leqslant n}|f_k(R_{\alpha,1})-R_{\alpha,k}|<\Phi(\beta_{\alpha})\}.$$
Denote by $\mathcal{R}_1$ the family of first coordinates of the points in $\mathcal{R_C}(\Phi)$; that is $$\mathcal{R}_1=(R_{\alpha,1})_{\alpha\in\mathcal{J_C}(\Phi)}.$$

\noindent\textsc{\textbf{Definition:} }\textit{The system $(\mathcal{R_C}(\Phi), \beta)$ is called locally ubiquitous with respect to $\rho$ if the system $(\mathcal{R}_1,\beta)$ is locally ubiquitous in $I$ with respect to $\rho$.}

Now, given an approximating function $\Psi$, let $\Lambda(\mathcal{R_C}(\Phi), \beta, \Psi)$ denote the set $x\in I$ for which the system of inequalities
\[\begin{cases}
|x-R_{\alpha,1}|<\Psi(\beta_{\alpha}) \\[2ex]
\max\limits_{2\leqslant k\leqslant n}|f_k(x)-R_{\alpha,k}|<\Psi(\beta_{\alpha})+\Phi(\beta_{\alpha}),
\end{cases}\]
is simultaneously satisfied for infinitely many $\alpha\in\mathcal{J}$. We will make use of the following lemmas.

\begin{lemma} \label{l01}
Consider the curve $\mathcal{C}=\mathcal{C}_f := \{(x, f_2(x), \ldots , f_n(x)) : x \in I\}$, where $f_2, \ldots, f_n$ are locally bi-Lipshitz in a finite interval $I$. Let $\Phi$ and $\Psi$ be approximating functions. Suppose that $(\mathcal{R}_{\mathcal{C}}(\Phi), \beta)$ is a locally ubiquitous system with respect to $\rho$. Let $s \in (0, 1]$ and suppose that $\Psi(2^{t+1}) \leqslant \frac{1}{2}\Psi(2^t)$ for $t$ sufficiently large. Then $$
\cH^s(\Lambda(\mathcal{R}_{\mathcal{C}}(\Phi), \beta,	\Psi)) = \cH^s(I) \;\quad {\rm if} \;\quad
\sum\limits_{t=1}^{\infty}\frac{(\Psi(2^t))^s}{\rho(2^t)}= \infty.$$

\end{lemma}

This is Lemma 4 in \cite{BDV_main}. For $s\neq 1$ the following is its generalization in terms of general Hausdorff measure and a consequence of Corollary 3 in \cite{BDV_mtl}.

\begin {lemma} \label{l02}
Consider the curve $\mathcal{C}=\mathcal{C}_f=\{(x,f(x)):x\in I\}$, $f\in C^{(3)}(I)$, locally Lipshitz. Let $\Phi$ and $\Psi$ be approximating functions and suppose that $(\mathcal{R_C}(\Phi),\beta)$ is locally ubiquitous with respect to $\rho$. Let $h$ be a regular dimension function such that $r^{-1}h(r)\to \infty$ as $r\to 0$ and $r^{-1}h(r)$ is decreasing. Suppose also $\Psi(2^{t+1}) \leqslant \frac{1}{2}\Psi(2^t)$ for $t$ sufficiently large. Then $$\cH^h(\Lambda(\mathcal{R_C}(\Phi),\beta,\Psi))=\infty\quad {\rm if}\quad \sum\limits_{n=1}^{\infty} h(\Psi(2^n))(\rho(2^n))^{-1}=\infty.$$
\end{lemma}

The following ubiquity statement is at the heart of establishing the divergence part of Theorem~\ref{t04}. It is a natural generalization of the homogeneous ubiquity statement established in \cite{BDV_main} (Corollary 5).

\begin{lemma} \label{l14}
Let $I$ denote a finite, open interval of $\R$ and let $f$ be a function in $C^{(3)}(I)$
such that $$
0 < c_1 := \inf_{x\in I} |f''(x)| \leqslant c_2 := \sup_{x\in I} |f''(x)| < \infty.$$
 Let $\psi$ be an approximating function satisfying
$$
\lim_{t\to +\infty}\psi(t) = \lim_{t\to +\infty}\frac{1}{t\psi(t)}= 0
$$
  and $\mathcal{C}=\mathcal{C}_f := \{(x, f(x)) : x \in I\}$ be a non-degenerate curve. Set
$$\beta : \mathcal{J} := \Z^2 \times \N \to \N : (\vv p, q) \to q ,\;\;\; \Phi : t \to t^{-1}\psi(t)\;\; and\;\; \rho : t \to \frac{u(t)}{t^2\psi(t)},$$
where $u : \R^{+} \to \R^{+}$ is any function such that $\lim_{t\to\infty}u(t) = \infty$. Given $\bm\theta\in\R^2$ denote by $\Q^2_{\mathcal{C}, \bm\theta}(\Phi)$ the set of resonant points from $\Q^2$ with numerators shifted by $\bm\theta$ and $\Phi$-close to the curve $\mathcal{C}$, i.e. $$\Q^2_{\mathcal{C}, \bm\theta}(\Phi)=\left\{\vv p/q\in\Q^2: \left|f\left(\frac{p_1+\theta_1}{q}\right)-\frac{p_2+\theta_2}{q}\right|<\Phi(\beta_{\alpha})\right\}.$$ Then the system $(\Q^2_{\mathcal{C},\bm\theta}(\Phi), \beta)$ is locally ubiquitous in $I^2$ with respect to $\rho$.
\end{lemma}

\noindent The key for establishing Lemma~\ref{l14} is the following statement (Theorem 4 from \cite{BVV}):

\begin{theoremBVV}\label{pt17}
Let $I$ denote a finite open interval of $\R$. Let $f$ be a function in $C^{(3)}(I)$ such that
$$0 < c_1\leqslant |f''(x)| \leqslant c_2 < \infty$$
 for some constants $c_1$ and $c_2$. Then for any interval $J\subseteq I$ there are constants $k_1, k_2, \delta_0, C_1, Q_0 > 0$ such that for any choice of $\psi$ and $Q > Q_0$ subject to $$\frac{k_1}{Q}\leqslant \psi(Q)\leqslant k_2$$
one has $$\left|\bigcup_{\vv p/q\in A_Q(J,\; \bm\theta)}\left(B\left(\frac{p_1+\theta_1}{q}, \frac{C_1}{Q^2 \psi(Q)}\right)\cap J\right)\right|\geqslant \frac{1}{2}|J|,$$
where
\begin {multline*}
A_Q(J, \bm\theta) := \{\vv p/q\in \Q^2: Q/2<q\leqslant Q, (p_1+\theta_1)/q\in J,\\
{}\left|f\left(\frac{p_1+\theta_1}{q}\right)-\frac{p_2+\theta_2}{q}\right|<\psi(Q)/Q\}.
\end {multline*}

\end{theoremBVV}

\begin{proof}[Proof of Lemma~\ref{l14}.]
Take a subinterval $J\subseteq I$ and define the set \begin {multline*}
A^{*}_Q(J, \bm\theta) := \{\vv p/q\in \Q^2: Q/u(Q)<q\leqslant Q, (p_1+\theta_1)/q\in J,\\
{}\left|f\left(\frac{p_1+\theta_1}{q}\right)-\frac{p_2+\theta_2}{q}\right|<\psi(Q)/Q\}.
\end {multline*}
Then, as $\lim_{t\to\infty}u(t) = \infty$, for $Q$ sufficiently large we have $1/u(Q)<1/2$, therefore $A_Q(J, \bm\theta)\subseteq A^{*}_Q(J, \bm\theta)$. Then,  applying the above theorem gives
\begin{multline*}\left|\bigcup_{\vv p/q\in A^{*}_Q(J,\; \bm\theta)}\left(B\left(\frac{p_1+\theta_1}{q}, \frac{u(Q)}{Q^2 \psi(Q)}\right)\cap J\right)\right|\geqslant \\
{} \geqslant\left|\bigcup_{\vv p/q\in A_Q(J,\; \bm\theta)}\left(B\left(\frac{p_1+\theta_1}{q}, \frac{C_1}{Q^2 \psi(Q)}\right)\cap J\right)\right|\geqslant \frac{1}{2}|J|,
\end{multline*}
which establishes the lemma.
\end{proof}

\section{Proof of Theorem~\ref{t04}}
The proof splits naturally into two parts: the convergence and divergence cases. We begin with the former.

\subsection{Convergence case.}

The proof of the convergence case follows on constructing a cover of the set $\mathcal{C}_f \cap \mathcal{A}(\psi_1, \psi_2, \bm\theta)$ with some bounded intervals, estimating the measure of each of them and their number and thereby finding an estimate for the Hausdorff measure of the entire set.

There is no loss of generality in assuming that $\psi_1(q)\leqslant\psi_2(q)$ for all $q$. Otherwise define $\psi(q)=\min\{\psi_1(q), \psi_2(q)\},\; \varphi(q)=\max\{\psi_1(q), \psi_2(q)\}$. It is clear that $$\mathcal{A}(\psi_1,\psi_2, \bm\theta)\subset \mathcal{A}(\psi,\varphi, \bm\theta)\cup \mathcal{A}(\varphi,\psi, \bm\theta),$$
therefore, to prove the convergence part of Theorem~\ref{t04} it is  sufficient to show that both sets $\mathcal{C}_f\cap\mathcal{A}(\psi,\varphi,\bm\theta)$ and $\mathcal{C}_f\cap\mathcal{A}(\varphi,\psi,\bm\theta)$ are of Hausdorff measure zero. We will consider one of these two sets - the other case is similar.
Thus, suppose $\psi_1(q)\leqslant\psi_2(q)$ and
$$
\sum = \sum\limits_q q\, h\left(\frac{\psi_1(q)}{q}\right)\psi_2(q)<\infty.$$

\noindent Let $\Omega_{f,\psi_1,\psi_2, \bm\theta}$ be the set of $x\in I$ such that the system of inequalities:
\begin{equation}\label{sysdef}
\begin {cases}
 |x-\frac{p_1+\theta_1}{q}|<\frac{\psi_1(q)}{q},\\[2ex]
 |f(x)-\frac{p_2+\theta_2}{q}|<\frac{\psi_2(q)}{q}
\end {cases}
\end{equation}
holds for infinitely many $q\in\N$ and $p_1,p_2\in\Z$.
Since the function $f$ is continuously differentiable, the map $x\mapsto (x,f(x))$ is locally bi-Lipshitz and so $$\cH^h(\mathcal{C}_f \cap \mathcal{A}(\psi_1, \psi_2, \bm\theta))=0 \quad\Longleftrightarrow \quad\cH^h(\Omega_{f,\psi_1,\psi_2, \bm\theta})=0.$$ Therefore it suffices to show that $\cH^h(\Omega_{f,\psi_1,\psi_2, \bm\theta})=0.$

Write $\sigma(\frac{\vv p}{q},\bm\theta)$ for the set of $x\in I$ satisfying \eqref{sysdef} for fixed $\frac{\vv p}{q}\in \Q^2$ and fixed $\bm\theta\in\R^2$. For $n\geqslant 0$ define
$$\Omega_{f,\psi_1,\psi_2, \bm\theta}(n)=\bigcup\limits_{\vv p/q\in \Q^2,\; \sigma(\vv p/q, \bm\theta)\neq\emptyset,\; 2^n\leqslant q<2^{n+1}} \sigma(\vv p/q,\bm\theta).$$
Then for each $l\in\N$ and $n\geqslant l$ the collection $\Omega_{f,\psi_1,\psi_2, \bm\theta}(n)$ is a cover of $\Omega_{f,\psi_1,\psi_2, \bm\theta}$ by the sets $\sigma(\vv p/q,\bm\theta)$.
Therefore, if we can estimate the size of each set $\sigma(\frac{\vv p}{q},\bm\theta)$ and the number of such sets, it will give us a bound for Hausdorff measure of $\mathcal{C}_f \cap \mathcal{A}(\psi_1, \psi_2, \bm\theta)$.

We may assume without loss of generality that $f$ is bounded on the interval $I$. Clearly $$\text{\upshape{diam}}\left(\sigma(\frac{\vv p}{q},\bm\theta)\right) \leqslant \frac{2\psi_1(q)}{q},$$
and we only need to find an upper bound for the number of the sets $\sigma(\frac{\vv p}{q},\bm\theta)$ in our cover. In other words, we need to estimate $$\#\{\vv p/q\in \Q^2: 2^n\leqslant q<2^{n+1}, \sigma(\vv p/q, \bm\theta)\neq\emptyset\},$$ where the symbol $\#$ denotes the cardinality of the set.

Assume $\sigma(\frac{\vv p}{q},\bm\theta)\neq \emptyset$ and let $x\in \sigma(\frac{\vv p}{q},\bm\theta)$.
We also require $f'$ to be bounded on $I$. This is true if $f''$ is bounded and $I$ is a bounded interval.

Let $B=\{x\in I: |f''(x)|=0\}$. By the conditions of the theorem we have $\cH^h(B)=0$.
There is no loss of generality in assuming that $I$ is an open interval in $\R$. The set $B$ is closed in $I$ and thus the set $G=I\setminus B$ is open. A standard argument then allows one to write $G$ as a countable union of bounded intervals $I_i$ on which
\begin{equation}
 0<c_1=\inf|f''(x)|\leqslant c_2=\sup|f''(x)|<\infty,
\label{1}
\end{equation}
 and assume that $|I_i|_\R \leqslant 1$. Since $f\in C^{(3)}(I)$, it follows that $|f''(x)-f''(y)|\ll|x-y|\leqslant|x-y|^\xi$ for any $x,y \in I_i$ and $0 \leqslant \xi \leqslant 1$; i.e. $f''\in Lip_\xi(I_i)$. In particular we may take $$1>\xi>\frac{3\eta-1}{1+\eta}.$$

Suppose that $\cH^h(\Omega_{f,\psi_1,\psi_2, \bm\theta} \cap I_i)=0$ for all $i\in \N$. Using the fact that $\cH^h(B)=0$ we have that
$$
\begin{array}{rcl}
\displaystyle \cH^h(\Omega_{f,\psi_1,\psi_2, \bm\theta}) & \leqslant & \cH^h(B \cup \bigcup\limits_{i=1}^{\infty}(\Omega_{f,\psi_1,\psi_2, \bm\theta} \cap I_i))\\[3ex]
& \leqslant & \displaystyle \cH^h(B) + \sum\limits_{i=1}^{\infty}\cH^h(\Omega_{f,\psi_1,\psi_2, \bm\theta} \cap I_i)=0.
\end{array}
$$
Thus, without loss of generality we can assume that $f$ satisfies \eqref{1} on $I$ and $I$ is bounded.

By the mean value theorem $f(x)=f\left(\frac{p_1+\theta_1}{q}\right)+f'(\widetilde{x})\left(x-\frac{p_1+\theta_1}{q}\right)$ for some $\widetilde{x}\in I$. Let $2^n\leqslant q < 2^{n+1}$ for now. By \eqref{sysdef} $$
\begin{array}{rcl}
\ds\left|f\left(\frac{p_1+\theta_1}{q}\right)-\frac{p_2+\theta_2}{q}\right| & \leqslant &\ds \left|f(x)-\frac{p_2+\theta_2}{q}\right|+\left|f'(\widetilde{x})\left(x-\frac{p_1+\theta_1}{q}\right)\right|\\[3ex]
& < & \displaystyle \frac{c_3\psi_2(q)}{q}\quad \leqslant\quad\frac{c_3\psi_2(2^n)}{2^n},
\end{array}
$$
 where $c_3$ is a positive constant.
Thus,
\begin {multline*}
\#\{\vv p/q\in \Q^2: 2^n\leqslant q<2^{n+1}, \sigma(\vv p/q, \bm\theta)\neq\emptyset\}\\
{}\leqslant \#\left\{\vv p/q\in \Q^2: 2^n\leqslant q<2^{n+1}, \frac{p_1+\theta_1}{q}\in I,\right.\qquad\qquad\qquad\qquad\qquad\;\;\\
{}\left.\qquad\qquad\qquad\qquad\qquad\qquad\qquad \left|f(\frac{p_1+\theta_1}{q})-\frac{p_2+\theta_2}{q}\right|<c_3\psi_2(2^n)/2^n\right\}\\
{}\quad\;\;\leqslant \#\,\{a/q\in \Q: q\leqslant 2^{n+1}, (a+\theta_1)/q\in I,\; ||q f(\frac{a+\theta_1}{q})-\theta_2||<2c_3\psi_2(2^n)\}.
\end {multline*}

\noindent The crucial step in estimating this expression is contained in the following lemma (Theorem 3 from \cite{BVV}):

\begin{lemma} \label{l04}
Suppose that $0 < \xi < 1$ and $f''\in Lip_{\xi}([\eta, \nu])$ and that $Q \geqslant 1$ and $0 < \delta < \frac{1}{2}$. Denote
\begin{multline*}
N(Q,\delta, \bm\theta)=
{} \#\,\{(a, q)\in \Z\times\N : q \leqslant Q, (a+\theta_1)/q\in I, ||qf((a+\theta_1)/q)-\theta_2||< \delta\}.
\end{multline*}
Then
$$ N(Q, \delta, \bm\theta) \ll \delta Q^2 + \delta^{-\frac{1}{2}}Q^{\frac{1}{2}+\varepsilon} + \delta^{\frac{\xi-1}{2}}Q^{\frac{3-\xi}{2}}.$$
\end{lemma}

In view of the choices of $\psi_2$ and $\xi$, Lemma~\ref{l04} implies that
\begin {equation}
\#\,\{\vv p/q\in \Q^2: 2^n\leqslant q<2^{n+1}, \sigma(\vv p/q, \bm\theta)\neq\emptyset\}\ll \psi_2(2^n)2^{2n}.
\end {equation}

Then, from the definition of Hausdorff measure
$$
\begin{array}{rcl}
\cH^h_l(\Omega_{f,\psi_1,\psi_2, \bm\theta}) & \leqslant & \displaystyle \sum\limits_{n=l}^{\infty}\quad \sum\limits_{\vv p/q\in \Q^2,\; \sigma(\vv p/q, \bm\theta)\neq\emptyset,\; 2^n\leqslant q<2^{n+1}} h(2\psi_1(2^n)/2^n) \\[4ex]
& \ll & \displaystyle \sum\limits_{n=l}^{\infty}2^{2n}h\left(\frac{\psi_1(2^n)}{2^n}\right)\psi_2(2^n)\to 0
\end{array}
$$
as $l\to \infty$. Hence $\cH^h(\Omega_{f,\psi_1,\psi_2, \bm\theta})=0$ and this completes the proof of the convergent case of the theorem.

\subsection{Divergent case.}
The divergent case follows relatively easily from the ubiquity lemmas (Lemma~\ref{l02} and Lemma~\ref{l14}). In short we need to show that all the conditions of Lemma~\ref{l02} are satisfied.

Again, there is no loss of generality in assuming that $\psi_1(q)\leqslant\psi_2(q)$ for all $q$.
To see this note that otherwise
$$
\begin{array}{rcl}
\ds\sum & := & \ds\sum\limits_q q\cdot h\left(\frac{\min\{\psi_1(q),\psi_2(q)\}}{q}\right)\max\{\psi_1(q),\psi_2(q)\}\\[2ex]
& = & \ds\sum\limits_{q: \psi_1(q)\leqslant \psi_2(q)} q\, h\left(\frac{\psi_1(q)}{q}\right)\psi_2(q)+\sum\limits_{q: \psi_1(q)>\psi_2(q)} q\, h\left(\frac{\psi_2(q)}{q}\right)\psi_1(q).
\end{array}
$$

Since the sum $\Sigma$ diverges, one of the components must diverge as well. Suppose that the first component diverges (for the second one the proof will be the same just with the indices interchanged). In other words,
$$\sum\limits_{q: \psi_1(q)\leqslant \psi_2(q)} q\cdot h\left(\frac{\psi_1(q)}{q}\right)\psi_2(q)=\infty.$$
This means that we can choose the infinite sequence $q_i$ from the set of $q: \psi_1(q)\leqslant \psi_2(q)$. Clearly, for all $i$ the condition $\psi_1(q_i)\leqslant\psi_2(q_i)$ will be satisfied.
Now construct the new function $\widetilde{\psi_1}(q)=\psi_1(q_{i+1})$ for $q_i<q\leqslant q_{i+1}$. Since $\psi_1$ is decreasing we can conclude that $\widetilde{\psi_1}(q)\leqslant \psi_1(q)$ and since the functions $\psi_1$ and $\psi_2$ are monotonic then $\widetilde{\psi_1}(q)\leqslant \psi_2(q)$ for all $q$. The set $\mathcal{A}(\widetilde{\psi_1}, \psi_2,\bm\theta)\subset \mathcal{A}(\psi_1, \psi_2,\bm\theta)$. This means that if $$\cH^h(\mathcal{A}(\widetilde{\psi_1}, \psi_2,\bm\theta))=\infty$$ then $$\cH^h(\mathcal{A}(\psi_1, \psi_2,\bm\theta))=\infty.$$
So, $\psi_1\leqslant\psi_2$ can be assumed.

Now since the sum $\Sigma$ diverges we can find a strictly increasing sequence of positive integers $\{l_i\}_{i\in\N}$ such that $$\sum\limits_{l_{i-1}<l\leqslant l_i} l\cdot h\left(\frac{\psi_1(l)}{l}\right) \psi_2(l)>1.$$
Define the function $u: l \to u(l):=i$ for $l_{i-1}<l\leqslant l_i$.

\noindent Note that
$$
\begin {array}{rcl}
\displaystyle\sum\limits_{l=1}^{\infty}l\cdot h\left(\frac{\psi_1(l)}{l}\right) \psi_2(l) (u(l))^{-1} & = &\displaystyle \sum\limits_{i=1}^{\infty} \sum\limits_{l_{i-1}<l\leqslant l_i}l\cdot h\left(\frac{\psi_1(l)}{l}\right) \psi_2(l) (u(l))^{-1}\\[4ex]
& > & \displaystyle\sum\limits_{i=1}^{\infty}i^{-1}\quad =\quad \infty
\end{array}
$$
and since $\frac{h(\psi_1) \psi_2}{u}$ is decreasing
$$
\begin {array}{rcl}
\infty & = & \ds\sum\limits_{t=0}^{\infty} \sum\limits_{2^t\leqslant l<2^{t+1}}l\cdot h\left(\frac{\psi_1(l)}{l}\right) \psi_2(l) (u(l))^{-1}\\
& \ll & \ds\sum\limits_{t=0}^{\infty}2^{2t}h\left(\frac{\psi_1(2^t)}{2^t}\right)\psi_2(2^t)(u(2^t))^{-1}.
\end {array}
$$
Hence,

\begin{equation}\label{lastone}
\sum\limits_{t=0}^{\infty}2^{2t}h\left(\frac{\psi_1(2^t)}{2^t}\right)\psi_2(2^t)(u(2^t))^{-1}=\infty.
\end{equation}

Now let $\Psi(t)=\frac{\psi_1(t)}{t}, \; \Phi(t)=\frac{\psi_2(t)}{t}$ and $\rho(t)=\frac{u(t)}{t^2\psi_2(t)}$.
By Lemma~\ref{l14}, the system $(\Q_{\mathcal{C},\bm\theta}^2(\Phi), \beta)$ is locally ubiquitous relative to $\rho$, where $\beta$ is given by
$$
\beta:\mathcal{J}:=\Z^2\times \N\to\N: (\vv p,q)\mapsto q.
$$

\noindent In view of \eqref{lastone}

$$\sum\limits_{t=0}^{\infty}\frac{h(\Psi(2^t))}{\rho(2^t)}= \sum\limits_{t=0}^{\infty}2^{2t}h\left(\frac{\psi_1(2^t)}{2^t}\right)\psi_2(2^t)(u(2^t))^{-1}=\infty$$ and since $\psi_1$ is decreasing $$
\Psi(2^{t+1})=\frac{\psi_1(2^{t+1})}{2^{t+1}}\leqslant \frac{1}{2}\cdot\frac{\psi_1(2^t)}{2^t}=\frac{1}{2}\Psi(2^t).
$$

Thus, the conditions of Lemma~\ref{l02} are satisfied and the set $\Lambda(\Q_{\mathcal{C},\bm\theta}^2(\Phi),\beta,\Psi)$  is of full measure. By definition, the latter set consists of points $x\in I$ such that
\[\begin {cases}
 |x-\frac{p_1+\theta_1}{q}|<\Psi(\beta_\alpha)=\frac{\psi_1(q)}{q}<\frac{2\psi_1(q)}{q},\\[2ex]
 |f(x)-\frac{p_2+\theta_2}{q}|<\Phi(\beta_\alpha) +\Psi(\beta_\alpha)=\frac{\psi_2(q)}{q}+\frac{\psi_1(q)}{q}<\frac{2 \psi_2(q)}{q}
\end {cases}\]
has infinitely many solutions $\vv p/q\in\Q^2$. Obviously, if $x\in \Lambda(\Q^2_{\mathcal{C},\bm\theta}(\Phi),\beta,\Psi)$ then the point $(x,f(x))$ is in $\mathcal{A}(2\psi_1,2\psi_2,\bm\theta)$ and to complete the proof we apply what has already been proven to the approximating functions $\frac12\psi_1$ and $\frac12\psi_2$.

\section{Proof of Theorem~\ref{t02}}
\subsection{Convergent case.}
For the sake of convenience, let $\psi=\psi_1$ and $\phi=\psi_2$. It is clear that $$ \mathcal{A}(\psi,\phi,\bm\theta)\subset \mathcal{A}(\psi^{*},\psi_{*},\bm\theta)\cup\mathcal{A}(\psi_{*},\psi^{*},\bm\theta),$$
where
$$\psi_{*}(q)=\min\{\psi(q),\phi(q)\}\quad\text{and}\quad\psi^{*}(q)=\max\{\psi(q),\phi(q)\}.$$

Since $\psi^{*}\psi_{*}=\psi\phi$, we have that $\sum\psi^{*}(q)\psi_{*}(q)<\infty$. Thus, to prove the convergence part of Theorem \ref{t02} it is  sufficient to show that both the sets $\mathcal{C}_f\cap\mathcal{A}(\psi^{*},\psi_{*},\bm\theta)$ and $\mathcal{C}_f\cap\mathcal{A}(\psi_{*},\psi^{*},\bm\theta)$ are of Lebesgue measure zero. We will consider one of these two sets - the other case is similar. Thus, without any loss of generality we assume that $\psi(q)\geqslant\phi(q)$ for all $q\in\N$.

Since $\sum\psi(q)\phi(q)<\infty$ and both $\psi$ and $\phi$ are decreasing, we have that $\psi(q)\phi(q)<q^{-1}$ for all sufficiently large $q$. Hence $\phi(q)\leqslant q^{-1/2}$ for sufficiently large $q$. Further, we can assume that \begin{equation}\label{02-1}
\psi(q)\geqslant q^{-2/3}
\end{equation}
for all $q\in\N$. To see this consider the auxiliary function $\widetilde{\psi}$ such that $$\widetilde{\psi}=\max\{\psi(q),q^{-2/3}\}.$$
Clearly, $\widetilde{\psi}$ is an approximating function. It also satisfies the following set inclusion:
$$\mathcal{A}(\psi,\phi,\bm\theta)\subset\mathcal{A}(\widetilde{\psi},\phi,\bm\theta).$$
Moreover,
$$
\begin{array}{rcl}
\displaystyle\sum\limits_{q=1}^{\infty}\widetilde{\psi}(q)\phi(q) & \leqslant & \displaystyle\sum\limits_{q=1}^{\infty}\psi(q)\phi(q) + \sum\limits_{q=1}^{\infty}q^{-2/3}\phi(q)\\ [4ex]
 & \ll & \displaystyle\sum\limits_{q=1}^{\infty}\psi(q)\phi(q) + \sum\limits_{q=1}^{\infty}q^{-2/3}q^{-1/2}\quad <\quad \infty.
\end{array}
$$

This means that it is sufficient to prove the statement with $\psi$ replaced  by $\widetilde{\psi}$ and therefore \eqref{02-1} can be assumed.

The set $\mathcal{C}_f \cap\mathcal{A}(\psi,\phi,\bm\theta)$ is a limsup set with the following natural representation: $$\mathcal{C}_f\cap\mathcal{A}(\psi,\phi,\bm\theta) = \bigcap\limits_{n=1}^{\infty} \bigcup\limits_{q=n}^{\infty} \bigcup_{(p_1,p_2)\in\Z^2} \mathcal{A}(p_1,p_2,q,\bm\theta),$$
where \begin{multline*}
\mathcal{A}(p_1,p_2,q,\bm\theta)=\left\{(x,y)\in\R^2: \left|x-\frac{p_1+\theta_1}{q}\right|<\frac{\psi(q)}{q},\right.\\
{}\left.\left|y-\frac{p_2+\theta_2}{q}\right|<\frac{\phi(q)}{q} \right\}.\qquad\qquad\qquad\;\;\;
\end{multline*}

Using the fact that $\psi$ is decreasing, it is readily verified for any $n\geqslant 1$ that
\begin{equation}\label{02-2}
\mathcal{C}_f\cap\mathcal{A}(\psi,\phi,\bm\theta)\subset \bigcup\limits_{t=n}^{\infty}\bigcup\limits_{2^t\leqslant q<2^{t+1}}\bigcup_{(p_1,p_2)\in\Z^2} \mathcal{C}_f\cap\mathcal{A}(p_1,p_2,q,\bm\theta,t),
\end{equation}
where
\begin{multline*}
\mathcal{A}(p_1,p_2,q,\bm\theta,t)=\left\{(x,y)\in\R^2: \left|x-\frac{p_1+\theta_1}{q}\right|<\frac{\psi(2^t)}{2^t},\right.\\
{}\left.\left|y-\frac{p_2+\theta_2}{q}\right|<\frac{\phi(2^t)}{2^t} \right\}\qquad\qquad\qquad\;\;\;
\end{multline*}
and $t$ is uniquely defined by $2^t\leqslant q<2^{t+1}$. As in the convergent part of the proof of Theorem~\ref{t04} we can assume without loss of generality that $c_1<f'(x)<c_2$ for all $x\in I$. It follows that
\begin{equation}\label{02-3}
|\mathcal{C}_f\cap\mathcal{A}(p_1,p_2,q,\bm\theta,t)|_{\mathcal{C}_f}\ll\frac{\phi(2^t)}{2^t}.
\end{equation}

Finally for a fixed $t$ let $N(t)$ denote the number of triples $(q,p_1,p_2)$ with $2^t\leqslant q<2^{t+1}$ such that $\mathcal{C}_f\cap\mathcal{A}(p_1,p_2,q,\bm\theta,t)\neq\emptyset$.
Suppose now that $\mathcal{C}_f\cap\mathcal{A}(p_1,p_2,q,\bm\theta, t)\neq\emptyset$. Then for some $(x,y)\in\mathcal{C}_f\cap\mathcal{A}(p_1,p_2,q,\bm\theta, t)$ and $\kappa_1,\kappa_2$ satisfying $-1<\kappa_1,\kappa_2<1$ we have
$$
x=\frac{p_1+\theta_1}{q}+\kappa_1 \frac{\psi(2^t)}{2^t},
$$
$$
y=\frac{p_2+\theta_2}{q}+\kappa_2 \frac{\phi(2^t)}{2^t}.
$$

Thus, it can be shown that $$
\begin{array}{rcl}
\ds f\left(\frac{p_1+\theta_1}{q}\right)-\frac{p_2+\theta_2}{q} & = & \ds f\left(\frac{p_1+\theta_1}{q}\right)-f(x)+f(x)-\frac{p_2+\theta_2}{q}\\[4ex]
 & = &\ds -\kappa_1 f'(\widetilde{x})\cdot \frac{\psi(2^t)}{2^t}+\kappa_2 \frac{\phi(2^t)}{2^t},
\end{array}
$$
where $\widetilde{x}$ lies between $x$ and $\frac{p_1+\theta_1}{q}$. Thus, one can deduce that
$$
\left|f\left(\frac{p_1+\theta_1}{q}\right)-\frac{p_2+\theta_2}{q}\right|\ll \frac{\psi(2^t)}{2^t} .
$$

\noindent On applying Lemma~\ref{l04} with $Q=2^{t+1}$, $\delta=\psi(2^t)$ and $\frac{3}{5}<\xi<1$ it follows that there exists an absolute constant $c>0$ such that

\begin{equation}\label{02-4}
N(t)\ll 2^{2t}\psi(2^t).
\end{equation}

The upshot of \eqref{02-2}, \eqref{02-3} and \eqref{02-4} is that
$$
\begin{array}{rcl}
|\mathcal{C}_f\cap\mathcal{A}(\psi,\phi,\bm\theta)|_{\mathcal{C}_f} & \ll & \displaystyle\sum\limits_{t=n}^{\infty}\sum\limits_{2^t\leqslant q<2^{t+1}}\sum_{(p_1,p_2)\in\Z^2}|\mathcal{C}_f\cap\mathcal{A}(p_1,p_2,q,\bm\theta)|_{\mathcal{C}_f}\\[4ex]
 & \ll & \displaystyle\sum\limits_{t=n}^{\infty}N(t)\frac{\phi(2^t)}{2^t}\\
  & \ll & \displaystyle\sum\limits_{t=n}^{\infty}2^t\psi(2^t)\phi(2^t)\\
 & \asymp & \displaystyle\sum\limits_{q=2^n}^{\infty}\psi(q)\phi(q).
\end{array}
$$

Since $\sum\limits_{q=1}^{\infty}\psi(q)\phi(q)<\infty$ we have that $\sum\limits_{q=2^n}^{\infty}\psi(q)\phi(q)\to 0$ as $n\to 0$ and it follows that
$$
|\mathcal{C}_f\cap\mathcal{A}(\psi,\phi,\bm\theta)|_{\mathcal{C}_f}=0
$$
as required.
This completes the proof of the convergent case of Theorem \ref{t02}.

\subsection{Divergent case.}
As $\mathcal{C}=\mathcal{C}_f$ is non-degenerate almost everywhere, we can restrict our attention to a sufficiently small patch of $\mathcal{C}_f$ which can be written as $\{(x,f(x)): x\in I_0\}$ where $I_0$ is a sub-interval of $I$ and
\begin{equation}\label{02d-1}
c_1<|f''(x)|<c_2\quad\text{for all } x\in I_0\text{ and some positive constants }c_1, c_2.
\end{equation}
Without loss of generality and for clarity we assume that $f$ satisfies \eqref{02d-1} on $I$.

We are given that $\psi_1$ and $\psi_2$ are approximating functions such that
\begin{equation}\label{02d-2}
\sum\limits_{q=1}^{\infty}\psi_1(q)\psi_2(q)=\infty.
\end{equation}
Thus, at least one of the following two sums diverges:
$$
\sum\limits_{q\in\N,\,\psi_1(q)\geqslant\psi_2(q)}\psi_1(q)\psi_2(q) \quad \text{or} \quad \sum\limits_{q\in\N,\,\psi_1(q)\leqslant\psi_2(q)}\psi_1(q)\psi_2(q).
$$
Throughout we assume that the sum on the right is divergent. The argument below can easily be modified to deal with the case that only the sum on the left is divergent.

Moreover, there is no loss of generality in assuming that
\begin{equation}\label{02d-3}
\psi_2(q)\geqslant\psi_1(q)\quad\text{for all }q\in\N.
\end{equation}
To see that, define the auxiliary function $\widetilde{\psi_1}:\; q\to \widetilde{\psi_1}:=\min\{\psi_1(q),\psi_2(q)\}.$ Then the sum
$$
\sum\limits_{q=1}^{\infty}\widetilde{\psi_1}(q)\psi_2(q)
$$
diverges since by assumption it contains a divergent sub-sum. It is readily verified that $\widetilde{\psi_1}$ is an approximating function and that $\mathcal{A}(\widetilde{\psi_1},\psi_2,\bm\theta)\subset\mathcal{A}(\psi_1,\psi_2,\bm\theta).$ Thus, to complete the proof of the divergent part of Theorem~\ref{t02} it suffices to prove the result with $\psi_1$ replaced by $\widetilde{\psi_1}$. Hence, without loss of generality, \eqref{02d-3} can be assumed.

Now there is no loss of generality in assuming that
\begin{equation}\label{02d-4}
\psi_i(q)\to 0\quad\text{as }\; q\to \infty\quad(i=1,2).
\end{equation}
To see that, define the increasing function $v:\R^{+}\to\R{+}$ as follows:
$$
v(q):=\sum\limits_{t=1}^{[q]}\psi_1(t)\psi_2(t).
$$
In view of \eqref{02d-2}, $\lim_{t\to\infty}v(t)=\infty$. Fix $k\in\N$, then
$$
\sum\limits_{t=k}^{m}\frac{\psi_1(t)\psi_2(t)}{v(t)}\geqslant \sum\limits_{t=k}^{m}\frac{\psi_1(t)\psi_2(t)}{v(m)}=\frac{v(m)-v(k-1)}{v(m)}\to 1\quad\text{as}\quad m\to\infty
$$
and so
$$
\sum\limits_{t=k}^{m}\frac{\psi_1(t)\psi_2(t)}{v(t)}\geqslant 1\quad\text{for all}\;\;k\;\;\text{and sufficiently large}\;\;m.
$$
This implies that the sum $\sum\limits_{t=1}^{\infty}\frac{\psi_1(t)\psi_2(t)}{v(t)}$ diverges. Next, for $i=1,2$ consider the functions
$$
\psi_i^{*}:q\to\psi_i^{*}(q):=\psi_i(q)/\sqrt{v(q)}.
$$
Then both $\psi_1^{*}(q)$ and $\psi_2^{*}(q)$ are decreasing and tending to $0$ as $q\to\infty$. Furthermore $\sum\limits_{q=1}^{\infty}\psi_1^{*}(q)\psi_2^{*}(q)=\infty$ and $\mathcal{A}(\psi_1^{*},\psi_2^{*},\bm\theta)\subset \mathcal{A}(\psi_1,\psi_2,\bm\theta)$. Therefore, it suffices to establish the statement for $\psi_1^{*},\psi_2^{*}$. Further, there is no loss of generality in assuming that
\begin{equation}\label{02d-5}
\psi_2(q)\geqslant q^{-2/3}\quad\text{for all } q.
\end{equation}
To this end, define $\widehat{\psi_2}(q)=\max\{\psi_2(q),q^{-2/3}\}$. In view of \eqref{02d-3} it is readily verified that
$$
\mathcal{A}(\psi_1,\widehat{\psi_2},\bm\theta)\subseteq \mathcal{A}(\psi_1,\psi_2,\bm\theta)\cup \mathcal{A}(q\mapsto q^{-2/3},q\mapsto q^{-2/3},\bm\theta).
$$
The set $\mathcal{A}(q\mapsto q^{-2/3},q\mapsto q^{-2/3},\bm\theta)$ is essentially the set $\mathcal{A}_{2}(q\mapsto q^{-2/3},\bm\theta)$. Therefore, the Khintchine-Jarn\'{i}k Theorem implies that the set $\mathcal{A}(q\mapsto q^{-2/3},q\mapsto q^{-2/3},\bm\theta)$ is of $2$-dimensional Lebesgue measure zero (which corresponds to $\mathcal{H}^2$ in the theorem) and hence
$$
|\{x\in I: (x,f(x))\in \mathcal{A}(\psi_1,\widehat{\psi_2},\bm\theta)\}|_{\mathcal{C}_f}\leqslant |\{x\in I: (x,f(x))\in \mathcal{A}(\psi_1,\psi_2,\bm\theta)\}|_{\mathcal{C}_f}.
$$
Thus, to complete the proof it suffices to prove that the set of the left has full measure. In turn, this justifies \eqref{02d-5}.

In view of \eqref{02d-4} and \eqref{02d-5} the function $\psi_2$ satisfies the condition
$$
\lim_{t\to+\infty}\psi_2(t)=\lim_{t\to+\infty}\frac{1}{t\psi_2(t)}=0
$$
and Lemma~\ref{l14} is applicable with $\psi=\psi_2$. By \eqref{02d-2} and the fact that $\psi_1$ and $\psi_2$ are decreasing we obtain that
$$
\infty=\sum\limits_{t=0}^{\infty}\sum\limits_{2^t\leqslant q<2^{t+1}}\psi_1(q)\psi_2(q) \leqslant \sum\limits_{t=0}^{\infty}\sum\limits_{2^t\leqslant q<2^{t+1}}\psi_1(2^t)\psi_2(2^t) = \sum\limits_{t=0}^{\infty}2^t\psi_1(2^t)\psi_2(2^t).
$$
Hence
\begin{equation}\label{02d-6}
\sum\limits_{t=0}^{\infty}2^t\psi_1(2^t)\psi_2(2^t)=\infty.
\end{equation}
Next, define the increasing function $u:\Rp\to\Rp$ as follows
$$
u(q)=\sum\limits_{t=0}^{[q]}2^t\psi_1(2^t)\psi_2(2^t).
$$
Trivially, $\lim_{t\to\infty}u(t)=\infty$. On using the same argument as when establishing \eqref{02d-4} above we verify that
\begin{equation}\label{02d-7}
\sum\limits_{t=0}^{\infty}\frac{2^t\psi_1(2^t)\psi_2(2^t)}{u(t)}=\infty.
\end{equation}

Now let $\Phi(t)=\frac{\psi_2(t)}{t}$ and $\rho(t)=\frac{u(t)}{t^2\psi_2(t)}$. By Lemma~\ref{l14}, $(\Q^2_{\mathcal{C},\bm\theta}(\Phi),\beta)$ is locally ubiquitous relative to $\rho$, where $\beta$ is given by
$$
\beta:\mathcal{J}:=\Z^2\times \N\to\N: (\vv p,q)\mapsto q.
$$
Let $\Psi(t)=\frac{\psi_1(t)}{t}$. In view of \eqref{02d-7}
$$
\sum\limits_{t=1}^{\infty}\frac{\Psi(2^t)}{\rho(2^t)}\quad = \quad \sum\limits_{t=1}^{\infty}\frac{\frac{\psi_1(2^t)}{2^t}}{\frac{u(2^t)}{2^{2t}\psi_2(2^t)}} \quad = \quad \sum\limits_{t=1}^{\infty}\frac{2^t\psi_1(2^t)\psi_2(2^t)}{u(t)} \quad = \quad \infty,
$$
and since $\psi_1$ is decreasing,
$$
\Psi(2^{t+1})=\frac{\psi_1(2^{t+1}}{2^{t+1}}\leqslant \frac{1}{2}\cdot\frac{\psi_1(2^t)}{2^t}=\frac{1}{2}\Psi(2^t).
$$
Thus, the conditions of Lemma~\ref{l01} are satisfied and it follows that the set $\Lambda(\Q^2_{\mathcal{C},\bm\theta}(\Phi),\beta,\Psi)$ is of full measure. By definition and \eqref{02d-3}, the latter set consists of points $x\in I$ such that the system
\[\begin {cases}
 |x-\frac{p_1+\theta_1}{q}|<\Psi(q)=\frac{\psi_1(q)}{q}<\frac{2\psi_1(q)}{q},\\[2ex]
 |f(x)-\frac{p_2+\theta_2}{q}|<\Phi(q) +\Psi(q)=\frac{\psi_2(q)}{q} + \frac{\psi_1(q)}{q}\leqslant\frac{2\psi_2(2)}{q}
\end {cases}\]
has infinitely many solutions $\vv p/q\in\Q^2$. Obviously, for $x\in \Lambda(\Q^2_{\mathcal{C},\bm\theta}(\Phi),\beta,\Psi)$ the point $(x,f(x))$ is in $\mathcal{A}(2\psi_1,2\psi_2,\bm\theta)$ and to complete the proof we simply apply what has already been proven to the approximating functions $\frac12\psi_1$ and $\frac12\psi_2$.

\section{Concluding Remarks}

Although we have considered a reasonably general approach in dealing with the problems in this paper, there are still some open problems which could improve the current state of knowledge concerning metric theory of inhomogeneous Diophantine approximation on planar curves. We summarize them here.

One possible improvement could be to remove the $C^3$ condition from Theorems \ref{t04} and \ref{t02} and replace it with $C^2$ or better, keeping in mind that the role of $C^3$ is mainly to keep the curve smooth enough to deviate from the hyperplane. In view of the perturbational approach used by authors in \cite{BZ} it might be possible to reduce the smoothness condition to $C^1$. In the case of divergence (both Theorems \ref{t04} and \ref{t02}) this can be done by using the ideas from paper \cite{BZ}. In the convergence case of Theorem \ref{t02} this can be done using the results from \cite{BVV}, but the convergence case of Theorem \ref{t04} will require the revision of results from \cite{BVV}, if at all possible.

\medskip

One can also consider a multiplicative varient of $\cA(\psi,\bm\theta)$; namely the set
$$\cA^{*}(\psi,\bm\theta):=\{(x_1,x_2)\in\R^2: \|qx_1-\theta_1\|\|qx_2-\theta_2\|<\psi(q)\text{ for infinitely many }q\in\N\}.$$
In the homogeneous setup, with $\bm\theta=(0,0)$, there is a criterion obtained by Gallagher \cite{Gallagher} concerning two-dimensional Lebesgue measure of this set.
\begin{theoremG}
Let $\psi:\N\to\R^+$ be a monotonic approximating function. Then
$$|\mathcal{A}^{*}(\psi)|=\begin {cases}
 0 \ &{\rm if } \quad \sum\limits_{q=1}^{\infty}(\log q)\,\psi(q)<\infty \\[3ex]
 \text{\upshape{Full}} \ &{\rm if } \quad \sum\limits_{q=1}^{\infty}(\log q)\,\psi(q)=\infty.
\end {cases}$$
\end{theoremG}

In \cite{BV}, Beresnevich and Velani determine a condition for the set $\mathcal{C}_f \cap A^*(\psi)$ to be of Lebesgue measure zero in the special case that $\mathcal{C}_f$ is a rational quadric $\mathcal{Q}$.  Subsequently,  Badziahin and Levesley (\cite{BL}) extended this to the case of arbitrary $C^{(3)}$ non-degenerate curves and to Hausdorff measures. Recently, Hussain and Yusupova (\cite{HTmulti}) generalized previously existing results to Inhomogeneous framework.
\begin{theorem}\label{t05}
Let $\psi$ be an approximating function, $\bm\theta=(\theta_1,\theta_2)\in\R^2$, $f\in C^{(3)}(I)$, where $I$ is an interval and let $\mathcal{C}_f=\{(x, f(x)): \; x\in I\}$ be a non-degenerate curve. Then for $\frac{2}{3}<s\leqslant 1$ we have that
$$\cH^s(\mathcal{C}_f\cap \mathcal{A}(\psi,\bm\theta))=0 \; \ \ \ {\rm if }\;\; \ \ \  \sum\limits_{q=1}^{\infty}q^{1-s}\psi^s(q)(\log q)^s<\infty.$$
\end{theorem}

Noticeably, all the results established so far in the multiplicative setup on planar curves are measure zero statements. An obvious next step is to prove  the divergent counterpart of Theorem \ref{t05}. That is, in general, one would like to show that

$$\cH^s(\mathcal{C}_f\cap \mathcal{A}(\psi,\bm\theta))=\infty \; \ \ \ {\rm if }\;\; \ \ \  \sum\limits_{q=1}^{\infty}q^{1-s}\psi^s(q)(\log q)^s=\infty.$$

\noindent We remark that by adapting the arguments in this paper and the ideas developed in \cite{BBV, BZ, HTmulti}, it is possible to prove the divergent part (stated above).  As a first step to solving the above problem Beresnevich and Velani in \cite{BVprep} has determined the Hausdorff dimension of the set $\mathcal{C}_f\cap \mathcal{A}(\psi,\bm\theta)$.

\medskip

\noindent \emph{\textbf{Acknowledgment.}} For their help, support and encouragement we would like to thank  Prof. Victor Beresnevich and Prof. Sanju Velani. We would like to thank the referee for the valuable comments.

\end {document}